\def\Bbb#1{{#1\kern-.800em #1}}
\def\Bbc#1{{#1\kern-.750em #1}}

\def \a{\alpha}
\def \b{\beta}
\def \g{\gamma}
\def \d{\delta}
\def \s{\sigma}
\def \eps{\epsilon}
\def \om{\omega}
\def \lam{\lambda}

\def \r{\right}
\def \tens{\mathop{\otimes}}

\def \>{\rangle}
\def \<{\langle}

\newtheorem{prop}{Proposition}
\newtheorem{lemma}{Lemma}

\def \be{\begin{equation}}
\def \eq{\end{equation}}
\def \bee{\begin{eqnarray}}
\def \eqq{\end{eqnarray}}
\def \been{\begin{eqnarray*}}
\def \eqqn{\end{eqnarray*}}
\def \n{\nonumber}

\def \R{{\Bbb R}}
\def \Z{{\Bbb Z}}
\def \C{{\Bbc C}}

\def \s {\;\;}
\def \A {{\cal A}}
\def \O {\Omega}
\def \R {{\cal R}}
\def \H {{\cal H}}
\def \xa {x^{a}}

\def \ka {\xi^{a}}
\def \kb {\xi^{b}}
\def \da {\chi_{a}}
\def \db {\chi_{b}}
\def \bd {\bar{\delta}}
\def \zb {\bar{z}}
\def \ab {\bar{a}}
\def \bb {\bar{b}}
\def \del {\partial}
\def \delb {\bar{\partial}}
\def \S {S_{q}^{2}}
\def \SU {SU_{q}(2)}
\def \T {\left(\begin{array}{ll}
                 \a & \b \\
                 \g & \d 
              \end{array}\right)}
\def \r {\rho}

\def \bo {\bar{0}}
\def \th {\theta}
\def \I {{\cal I}}

\def \eb {\bar{e}}

\def \Dir {{\cal D}}

\documentstyle[12pt]{article}
\input math_macros.tex
\double

\begin{document}
\begin{titlepage}
\begin{center}
Jan. 9, 96
          \hfill  LBL-37902\\
          \hfill  UCB-PTH-95/37\\
          \hfill  q-alg/9505021\\
\vskip .3in

{\large \bf Riemannian Geometry on Quantum Spaces}

\vskip .5in

{Pei-Ming Ho\footnote{email address:
pmho@physics.berkeley.edu}}

\vskip .2in
{\em
Department of Physics \\
University of California \\ 
and\\
Theoretical Physics Group\\
Lawrence Berkeley Laboratory\\
University of California\\
Berkeley, CA 94720}
\end{center}

\begin{abstract}

An algebraic formulation of Riemannian geometry on quantum spaces
is presented,
where Riemannian metric, distance, Laplacian, connection,
and curvature have their counterparts.
This description is also extended to complex manifolds.
Examples include the quantum sphere,
the complex quantum projective spaces
and the two-sheeted space.

\end{abstract}

\end{titlepage}

\section{ Introduction }

In \cite{CFF} Chamseddine, Felder and Fr\"{o}hlich developed
the notions of Riemannian metric, connection and curvature
in the framework of the non-commutative geometry of A. Connes \cite{Con}.
In their formulation the Hilbert space representation is a prerequisite.
The purpose of this paper is to propose
a purely algebraic formulation of Riemannian geometry
on quantum spaces.
It is suitable for physicists to build physical models.
The question of mathematical rigor is left for future study.

In Sec.\ref{Rieman} we describe this algebraic formulation of
Riemannian geometry on quantum spaces and quantum complex manifolds.
It is applied to the quantum sphere $S_{q}^{2}$\cite{CHZ}
in Sec.\ref{q-sphere}.
In particular the explicit expression for the quantum distance is worked out.
A comment on the connection with A. Connes' work is made.
The complex projective spaces $CP_{q}(N)$\cite{CHZ3}
is considered in Sec.\ref{CPN},
and the two-sheeted space \cite{CFF} in Sec.\ref{Sheets}.

\section{ Riemannian Structure on Quantum Spaces }\label{Rieman}

The general notion of non-commutative differential calculus
is reviewed in Sec.\ref{DC}.
Those who are familiar with it can start with Sec.\ref{RVT}.

\subsection{ Differential Calculus on Quantum Spaces }\label{DC}

A quantum space is specified by an unital, associative,
non-commutative $*$-algebra $\A$ 
generated by $\{1, \xa, \I_{i}\}$ over the field $k = \C$,
where $\xa$'s are the coordinates on the quantum space,
and $\I_{i}$'s are non-commutative constants,
including for example the generators of the algebra of functions
of the quantum group
which specifies the quantum symmetry on the quantum space.
$\A$ is called the algebra of functions on the quantum space.

To talk about differential geometry on the quantum space,
we should have $\A$ extended to the algerba of differential calculus $\O(\A)$ generated by $\{1, \xa, \ka, \da, \I_{i}\}$,
where the commutation relations among the generators are given
so that one knows how to rewrite a product of elements in $\O(\A)$
in any prefered order of the elements.
(But the commutation relations between $\ka$ and $\db$ are not necessary.)
The $\ka$'s are differential one-forms and the $\da$'s are 
the derivations dual to them
so that the exterior derivative is $d = \ka \da$.
The $*$-involution on functions is also extended to all elements in $\O(\A)$
and it always reverses the ordering of a product.
All constants, namely the unity or $\I_{i}$'s, should vanish under $d$.

All the commutation relations should never mix terms of differential forms 
of different degrees so that $\O(\A)$ is graded.
The action of $d$ is: $(da) := [d, a] = da - ad$
for forms of even degrees (including elements in $\A$),
and $(da) := \{d, a\} = da + ad$ for forms of odd degrees.
Leibniz rule follows this definition.
We take the convention that
$(da)^{*} = d(a^{*})$ for even degrees and $(da)^{*} = -d(a^{*})$ 
for odd degrees.
That is, $d$ is anti-self-dual: $d^{*} = -d$.
We also require the nilpotency of $d$, namely, $dd\equiv 0$.

\subsection{ Riemannian Metric, Vector Fields and Tensor Fields }
\label{RVT}

A general coordinate transformation is specified by:
$\xa \rightarrow x'^{a}$,
where $x'^{a} = x'^{a}(x)$ are elements in $\A$.
(Einstein's summation convention applies to the whole paper
unless otherwise stated.)
This transformation induces the transformation of $\ka$.
For example, if 
\be
   \ka = \sum f^{a}_{i}(x)dg^{a}_{i}(x), \label{ka}
\eq
then 
\[
   \ka \rightarrow \xi'^{a} = \sum f^{a}_{i}(x')dg^{a}_{i}(x'),
\]
where `$a$' is not summed over.
Re-expressing $\xi'^{a}$ in terms of $\ka$ and using commutation relations 
between $\xa$ and $\ka$ one can re-write the formula above as:
\[
   \ka \rightarrow \xi'^{a} = \kb M_{b}^{\s a}(x)
\]
for a certain matrix $M_{b}^{\s a}$ of elements in $\A$.

Since the transformation is not supposed to change the exterior derivative,
$d = \ka \da \rightarrow d' = \xi'^{a} \chi'_{a} = d$, 
so $\da \rightarrow {M^{-1}}_{a}^{\s b} \chi_{b}$.

A Riemannian metric $g^{ab}(x)$ is an invertible matrix of elements in $\A$ 
which transforms like a rank-two tensor (to be defined later):
\[
   g^{ab} \rightarrow g'^{ab} = M_{c}^{*a} g^{cd} M_{d}^{\s b},
\]
and is also Hermitian-symmetric:
\be
   (g^{ab})^{*} = g^{ba}. \label{gab}
\eq
Note that this symmetry is preserved by the transformation.
In the classical case there is no need of the $*$-involution
(complex conjugation) in (\ref{gab}) if all coordinates are real.
But if one is allowed to use complex coordinates,
for example, $(x, y)\rightarrow (x+iy, x-iy)$,
then (\ref{gab}) is a reality condition
for the Riemannian manifold.
The existence of the inverse of $g^{ab}$ is assumed
and it is denoted as $g_{ab}$ 
so that: 
$g^{ab} g_{bc} = \delta^{a}_{c} = g_{cb} g^{ba}$.
The transformation of $g_{ab}$ follows this definition.

A covariant vector field is a set of elements $\{\a_{a}\}$ in $\A$ 
which transform like 
$\{\da\}$:
\[
   \a_{a} \rightarrow M_{a}^{-1b} \a_{b}.
\]
Similarly, $\{\b^{a}\}$ is called a contra-variant vector field 
if it transforms like
$\{\ka\}$:
\[
   \b^{a} \rightarrow \b^{b} M_{b}^{\s a}.
\]
Note that $\a = \ka \a_{a}$ and $\b = \b^{a} \da$ are both invariant:
$\a \rightarrow \a$, $\b \rightarrow \b$,
so that we can simply use $\a, \b$ to denote the vector fields $\{\a_{a}\}$ 
and $\{\b^{a}\}$ in a coordinate-independent way.

Similarly we can define rank-two tensors of different types according to 
their transformations:
\been
   &\a^{ab} \rightarrow M_{c}^{*a} \a^{cd} M_{d}^{\s b}, \\
   &\a_{a}^{\s b} \rightarrow M_{a}^{-1c} \a_{c}^{\s d} M_{d}^{\s b}, \\
   &\a^{a}_{\s b} \rightarrow M_{c}^{*a} \a^{c}_{\s d} (M_{b}^{-1d})^{*}, \\
   &\a_{ab} \rightarrow M_{a}^{-1c} \a_{cd} (M_{b}^{-1d})^{*}.
\eqqn
Just like in the commutative case, the positions of indices of a tensor 
tell you the way it transforms.
In these formulas,
the ordinary contraction of indices is of this type: $\searrow$,
When indices are contracted as $\nearrow$, $\ast$-involution is involved.

Furthermore the Riemannian metric $g^{ab}$ can be used to raise or 
lower indices:
\bee
   &\a^{a} = \a_{b}^{*} g^{ba}, \label{raise} \\
   &\a_{a} = g_{ab} (\a^{b})^{*}, \label{lower} \\
   &\a^{ab} = g^{ac}\a_{c}^{\s b}, \n \\
   &\a_{ab} = \a_{a}^{\s c}g_{cb}, \n
\eqq
and so on.
Because of eq.(\ref{gab}),
if we raise and then lower an index we will get back to the original 
object.
It can also be checked that $(\a^a)^{*}\b^b$ and $\a_a\b^b$
are tensors if $\a$ and $\b$ are vectors.
The contraction of indices of two tensors can make a new tensor.
But sometimes the contraction of indices has to be accompanied by 
appropriate $*$-involution.

As $\a^{a}\b_{a}$ is invariant for any two vector fields $\a$ and $\b$,
we can always use the Riemannian metric to define the inner product 
$\<\cdot,\cdot\>$ between vector fields $\a$, $\b$ as in the classical case.
For example,
\bee
   &\<\a, \b\> = \<\ka\a_{a}, \kb\b_{b}\> 
            = \a_{a}^{*}\<\ka, \kb\>\b_{b}
            = \a_{a}^{*}g^{ab}\b_{b}, \n \\
   &\<\a, \b\> = \<\a^{a}\da, \b^{b}\db\> 
            = \a^{a}\<\da, \db\>\b^{b*} 
            = \a^{a}g_{ab}\b^{b*}. \n
\eqq
In both cases, $\<\a, \b\> = \<\b, \a\>^{*}$.

The magnitude of a vector $|\a|^{2} := \<\a, \a\>$ is real: $|\a|^{2*} = |\a|^{2}$,
due to eq.(\ref{gab}).

The invariant operator
\be
   \nabla^{2} := \chi^{a} \da = \chi_{a}^{*} g^{ab} \db \label{nabla}
\eq
is called the Laplacian.
It can be used to define the equation of motion for a scalar field $\Phi$ 
with mass $m$:
\[
   (\nabla^{2}+m^{2})\Phi = 0.
\]

{\em The non-commutativity forbids any tensor of rank higher than two.}
Therefore physical laws, if written as equations of motion, 
can only be written in terms of scalars, vectors, and rank-two tensors.
Fortunately all major classical physical laws are governed by
tensor equations of rank less than or equal to two.

Given a Hilbert space representation of the algerba $\A$ on $\H$,
one can define the ``distance'' between two states $s, s'$ (as generalized points) 
as \cite{Con}:
\be
   D(s, s') := sup\{|s(f)-s'(f)|; \quad \left\| |df|^{2} \right\| \leq 1, f\in \A\}. 
\label{dis}
\eq
The definition of $|df|^{2}$ is based on the metric $g^{ab}$ and therefore
the metric possesses the classical geometrical meaning.
The norm $\|\cdot\|$ is defined by: 
\[
   \|f\|^{2} := sup\{\frac{\<\phi|f^{*}f|\phi\>}{\<\phi|\phi\>};
                            \quad |\phi\>\in \H\}.
\]
An algebraic version of quantum distance can also be given
without mentioning Hilbert space representations.
An example is given in Sec.\ref{q-sphere}.

\subsection{ Connection }

As it is defined in \cite{CFF},
a connection $\nabla$ acts on $f \in \A$ as $d$:
$\nabla f = df$,
on $\ka$ by the connection one-forms $\om_{b}^{\s a}$:
\be
   \nabla \ka := \kb \otimes_{\A} \om_{b}^{\s a}, \label{xi-om}
\eq
and on a one-form $\a = \ka \a_{a}$ by Leibniz rule:
\been
   \nabla \a &=& (\nabla \ka)\a_{a} - \ka\otimes_{\A}(d\a_{a}) \\
             &=& -\ka \otimes_{\A} (\nabla \a)_{a},
\eqqn   
where $(\nabla \a)_{a} := d\a_{a} - \om_{a}^{\s b}\a_{b}$.

For $(\nabla\a)_{a}$ to be a covariant vector 
the connection one-form has the transformation:
\be
   \om_{a}^{\s b} \rightarrow M_{a}^{-1c}\om_{c}^{\s d}M_{d}^{\s b} 
   - M_{a}^{-1c}dM_{c}^{\s b}.
   \label{om}
\eq

For the Leibniz rule to hold on the inner product:
\[
   d\<\a, \b\> = \nabla\<\a, \b\> = (\nabla\a)^{a}\b_{a}+\a^{a}(\nabla\b)_{a},
\]
we define:
\[
   (\nabla \a)^{a} := d\a^{a} + \a^{b}\om_{b}^{\s a},
\]
which also ensures that $(\nabla \a)^{a}$ is a contra-variant vector.

The covariant derivation of $\a$ in the direction of the vector field $\b$ is:
\[
   (\nabla_{\b}\a)^{a} = \<\b, (\nabla\a)^{a}\>.
\]
The equation of geodesic flows is therefore:
\[
   (\nabla_{\a}\a)^{a} = 0.
\]

To define the action of $\nabla$ on rank-two tensor fields
we consider the scalar $f=\a_{a}^{*}\g^{ab}\b_{b}$.
Because $f$ is a scalar field, we have the equation: $\nabla f = df$.
By the undeformed Leibnitz rules of $\nabla$, 
we should have for the left hand side
(omitting the symbols $\otimes_{\A}$):
\been
   \nabla(\a_{a}^{*}\g^{ab}\b_{b}) &=& 
   (\nabla\a)_{a}^{*}\g^{ab}\b_{b}+\a_{a}^{*}(\nabla\g)^{ab}\b_{b}
   +\a_{a}^{*}\g^{ab}(\nabla\b)_{b} \\
   &=& (d\a_{a}-\om_{a}^{\s c}\a_{c})^{*}\g^{ab}\b_{b}
   +\a_{a}^{*}(\nabla\g)^{ab}\b_{b}
   +\a_{a}^{*}\g^{ab}(d\b_{b}-\om_{b}^{\s c}\b_{c}),
\eqqn
and for the right hand side:
\[
   d(\a_{a}^{*}\g^{ab}\b_{b}) = d(\a_{a}^{*})\g^{ab}\b_{b}
   +\a_{a}^{*}d\g^{ab}\b_{b}+\a_{a}^{*}\g^{ab}d\b_{b}.
\]
Identifying them we find:
\[
   (\nabla \g)^{ab} := d\g^{ab} + \g^{ac}\om_{c}^{\s b} 
   + (\om_{c}^{\s a})^{*}\g^{cb},
\]
which also ensures $(\nabla\g)^{ab}$ to be a rank-two tensor.

Suppose one has the physical law $(\nabla \a)^{a} = \b^{a}$,
it is equivalent to $(\nabla \a)_{a} = \b_{a}$ if
\be
   (\nabla g)^{ab} = dg^{ab} + \om^{ab} + (\om^{ba})^* = 0, \label{nablag}
\eq
which is called the metricity condition.

If $dg^{ab} = 0$, we have from the metricity condition $(\om^{ab})^{*} 
= -\om^{ba}$,
where $\om^{ab} := g^{ac}\om_{c}^{\s b}$.

The torsion $T^{a}$ is the covariant vector defined by \cite{CFF}:
\[
   T^{a} := (d - m\circ\nabla)\ka = d\ka - \kb\om_{b}^{\s a},
\]
where $m$ is the multiplication map $m(\a\otimes \b) := \a\b$.

In the classical case eq.(\ref{nablag}) and 
\be
   T^{a} = 0, \label{tor}
\eq
plus the reality conditions
imply that the connection one-form $\om_{a}^{\s b}$
is uniquely fixed by the metric $g^{ab}$.
The general expression for the analogous reality conditions
in the quantum case is so far unknown.
The difficulty is that the general transformation
will spoil the reality of a non-invariant quantity.
One has to invent appropriate conditions
for each particular case according to its algebraic properties.

For the quantum complex Hermitian manifolds defined later,
the situation is much simpler.
Just like their classical counterparts, 
Eq.(\ref{nablag}) alone determines the connection uniquely.

\subsection{ Curvature }

The curvature two-form is a rank-two tensor defined by:
\be
   R_{a}^{\s b} := d\om_{a}^{\s b} - \om_{a}^{\s c}\om_{c}^{\s b}. \label{CTF}
\eq
Using (\ref{nablag}), one can show that
\[ (R_{ab})^* = R_{ba}. \]

Is is easy to check that the Bianchi identity and consistency condition 
are satisfied:
\bee
   &dR_{a}^{\s b} - \om_{a}^{\s c}R_{c}^{\s b} 
   + R_{a}^{\s c}\om_{c}^{\s b} = 0, \n \\
   &dT^{a} = \xi^{b}R_{b}^{\s a}. \n
\eqq

Classically, in order to have the scalar curvature and Ricci tensor
one usually just strips the differential forms from 
the curvature two-form $R_{a}^{\s b}$ to get $R_{a\s cd}^{\s b}$
and then contracts $b, d$ for Ricci tensor,
and contracts in addition $a, c$ for the scalar curvature.
However in the deformed case this kind of operation is not 
covariant under general transformations.

Another more elegant way of defining the classical scalar curvature
and the Ricci tensor is to use the Hodge-$*$.
In the quantum case,
The Hodge-$*$ is required to satisfy
\be
   *(f\a g) = f(*\a)g \quad \forall f, g \in \A, \a \in \O(\A) \label{hs} 
\eq
and
\be
   (*\a)^* = *(\a^*) \quad \forall \a \in \O(\A), \label{starstar}
\eq
so that the scalar curvature defined as
\be
   \R := (-1)^{D+1} *(\xi^{a}(*R_{a}^{\s b})\xi_{b}) \label{hSC}
\eq
is invariant under general transformations
($D$ is the dimension of the space) and
real ($\R^* = \R$).
The integral
\be
\label{Lag}
\int \xi^a(*R_{a}^{\s b})\xi_{b}
\eq
is a candidate
for the action of a gravitational theory on the quantum space.

Similarly one can try to define the Ricci tensor as
\[
   \R_{a}^{\s b} := *((*R_{a}^{\s c})\xi_{c}\xi^{b}). \label{hRT}
\]
There are, however, many other inequivalent expressions
that are covariant under the general transformations.
For example, it is equally justified to define the Ricci tensor as
\[
   \R_{a}^{\s b} := *(\xi_{a}\xi^{c}(*R_{c}^{\s b})).
\]
This ambiguity in the Ricci tensor makes the scalar curvature better
for physical applications.

One can define the operator $\d := -*d*$ for a quantum space,
and naturally one will define the Laplacian by $\nabla^{2} := -(d+\d)^{2}$,
which is equivalent to $-\d d$ when acting only on functions if $\d\d = 0$.
In such cases the metric is determined by the Hodge-$*$ according to (\ref{nabla}).

\subsection{ Complex Manifolds } \label{Complex}

We define quantum complex manifolds
(more precisely, Hermitian manifolds) to be an associative $\ast$-algebra $\A$ 
generated by $\{1, z^{a}, \zb^{\ab}\}$ together with its differential calculus 
$\Omega(\A)$ generated by 
$\{1, z^{a}, \zb^{\ab}, dz^{a}, d\zb^{\ab}, \del_{a}, \delb_{\ab}\}$ 
with the following properties:
\begin{enumerate}
   \item $d = \d +\bd$, where $\d=dz^{a}\del_{a}$ and
         $\bd=\delb_{\ab}d\zb^{\ab}$ 
         with
         $\d\d = \bd\bd = 0$ and $\d\bd=-\bd\d$.
         $\d$ and $\bd$ should observe Leibniz rule separately,
         and $(\d\a)^{*} = (-1)^{p}\bd(\a^{*})$ for any form $\a$ of degree $p$.
   \item The generators of the algebra are divided into the holomorphic part
         $\{z^{a}\}$ and the anti-holomorphic part $\{\zb^{\ab}=(z^{a})^{*}\}$.
         The one-forms $\{dz^{a}\}$ and $\{d\zb^{\ab}\}$ are all independent.
   \item Denote $\Omega^{+}(\A):=\{dz^{a}\a_{a}: \a_{a}\in\A\}$ and
         $\Omega^{-}(\A):=\{\a_{\ab}d\zb^{\ab}: \a_{\ab}\in\A\}$.
         The commutation relations between one-forms and functions in $\A$
         are such that $\Omega^{\pm}(\A)\A = \A\Omega^{\pm}(\A)$.
         That is, $\Omega^{\pm}(\A)$ do not get mixed by commutation. 
   \item The metric $g^{\ab b}=(g^{\bb a})^{*}$ is given.
         The connection one-forms $\om_{a}^{\s b}$ and $\om_{\ab}^{\s \bb}
         :=(\om_{a}^{\s b})^{*}$ are also given
         such that $\om_{a}^{\s b}\in 
         \Omega^{+}(\A)$ and $\om_{\ab}^{\s \bb}\in \Omega^{-}(\A)$.
         Leibniz rule holds for the connection $\nabla$ and we have:
         $\nabla dz^{a} = dz^{b} \otimes_{\A} \om_{b}^{\s a}$ and
         $\nabla d\zb^{\ab} = -d\zb^{\bb} \otimes_{\A} \om_{\bb}^{\s \ab}$. 
\end{enumerate}

The coordinate transformations are restricted to the
holomorphic transformations only.
Holomorphic transformations are defined as those
which maps $z^{a}$ to holomorphic functions $f^{a}(z)$ and
$\zb^{\ab}=(z^{a})^{*}$ to $(f^{a}(z))^{*}$.
The properties of a complex manifold imply that this transformation
induces the map $dz^{a} \rightarrow dz^{b} M_{b}^{\s a}$
where $M_{a}^{\s b}$ is holomorphic, namely, $\bd M_{a}^{\s b}=0$.
Similarly, we have $d\zb^{\ab} \rightarrow (dz^{b} M_{b}^{\s a})^{*} = 
(M_{a}^{\s b})^{*}d\zb^{\bb}$
where $(M_{a}^{\s b})^{*}$ is anti-holomorphic, 
$\d (M_{a}^{\s b})^{*}=0$.

All the formulas we had before for Riemannian manifolds
can be easily modified for a complex manifold
with the understanding that the indices are only summed over the holomorphic 
or anti-holomorphic part.

With the fourth property of a complex manifold eq.(\ref{nablag}) says:
\[
   dg^{\ab b}+\om_{\bar{c}}^{\s \ab}g^{\bar{c}b}+g^{\ab c}\om_{c}^{\s b} = 0.
\]   
Since $d = \d+\bd$ we can separate the equation into $\Omega^{+}(\A)$ and 
$\Omega^{-}(\A)$,
and so the connection can be directly solved:
\be
   \om_{a}^{\s b} = -g_{a\bar{c}}(\d g^{\bar{c}b}). \label{omg1}
\eq
Only the holomorphic transformations will be consistent with this solution.
(That is, the transformation of the connection induced from 
the transformation of the metric by these expressions will be the same as
(\ref{om}) only for holomorphic transformations.)

From eq.(\ref{CTF}) the curvature two-form can now be expressed directly 
in terms of the metric as:
\be
   R_{a}^{\s b} = \bd\om_{a}^{\s b} = -\bd(g_{a\bar{c}}\d g^{\bar{c}b}).
   \label{Rg1}
\eq

The curvature two-form gives the scalar curvature according to (\ref{hSC})
with the indices restricted to the holomorphic or anti-holomorphic part,
and so it is in general not invariant
under general coordinate transformations.
The scalar curvature for a quantum complex manifold is
\be
   \R = (-1)^{D+1} *(dz^{a} (*R_{a\bb}) d\zb^{\bb}).
\eq
The definition of the Ricci tensor is not unique.

The condition (\ref{hs}) for the Hodge $*$ can be weakened to be no more than
an ordering prescription:
\[
   *(f(z)\a g(\zb)) = f(z)(*\a)g(\zb).
\]
Due to the extensive use of the $*$-involution in the quantum case,
we see that the complex structure helps to admit a Riemannian structure.

\section{ The Quantum Sphere $S_{q}^{2}$ }\label{q-sphere}

The Riemannian structure on a quantum space
with a quantum group symmetry should
respect its quantum symmetry.

In this section we describle one perticular quantum sphere $S_{q}^{2}$
in the family of quantum spheres of Podle\'{s} \cite{Pod}
with the $SU_{q}(2)$ symmetry (the one with $c=0$ in his notation)
in terms of the stereographic projection coordinates \cite{CHZ}
as an example of both Riemannian manifolds and complex manifolds.
For the convenience of the reader we review the stereographic 
projection of $S_{q}^{2}$ \cite{CHZ} in the following section.

\subsection{ The Stereographic Projection of $S_{q}^{2}$ }

$\S$ is the homogeneous space of $\SU/U(1)$ for $0 \leq q \leq 1$.
The algebra and differential calculus of $\S$
in terms of the complex coordinates $(z, \zb)$
can be induced from those of $\SU$ by the identification:
\bee
   z := \a\g^{-1}, & \zb := -\d\b^{-1}, \label{S-SU}
\eqq
where $\a,\b,\g,\d$ are the elements of an $\SU$-matrix $\T$.
Classically $z$, $\zb$ are the stereographic projection coordinates
projected from the north pole onto the tangent plane at the south pole.

The commutation relation between $z$, $\zb$ is therefore:
\[
   q(1+z\zb) = q^{-1}(1+\zb z),
\]
and the differential calculus induced from
the left-covariant 3D calculus on $\SU$ \cite{Wor1},
which is equivalent to one of the possible differential structures
on $\S$ studied by Podle\'{s} \cite{Pod3},
is specified by the following list of commutation relations:
\bee
   &z dz=q^{-2}dz z,     & \zb dz=q^{2}dz \zb,  \n \\
   &z d\zb=q^{-2}d\zb z, & \zb d\zb=q^{2}d\zb \zb, \n \\ 
   &dz dz=d\zb d\zb=0,  &dz d\zb=-q^{-2}d\zb dz, 
   \label{dzdzb} \\
   &\del z=1+q^{-2}z \del, &\del \zb=q^{2}\zb \del, \n \\
   &\delb z=q^{-2}z \delb, &\delb \zb=1+q^{2}\zb \delb, \n \\
   &\del \delb=q^{-2}\delb \del, \n 
\eqq
where $d=dz\del+d\zb\delb$ when acting on functions

The $*$-involution
on this whole algebra also follows from that of $\SU$:
\been
   &z^{*}=\zb, \\
   &(dz)^{*}=d\zb, \\ 
   &\del^{*}=-q^{-2}(1+\zb z)^{2}\delb(1+\zb z)^{-2}, \\
   &\delb^{*}=-q^{2}(1+\zb z)^{2}\del(1+\zb z)^{-2}.
\eqqn

The left $\SU$ transformation on $\SU$ induces an $\SU$ transformation on $\S$:
\been
   &z \rightarrow (az+b)(cz+d)^{-1}, \\
   &\zb \rightarrow (-c+d\zb)(a-b\zb)^{-1}, \\
   &dz \rightarrow dz(q^{-1}cz+d)^{-1}(cz+d)^{-1}, \\
   &d\zb \rightarrow d\zb(a-qb\zb)^{-1}(a-b\zb)^{-1}, \\
   &\del \rightarrow (cz+d)(q^{-1}cz+d)\del, \\
   &\delb \rightarrow (a-b\zb)(a-qb\zb)\delb,
\eqqn
where $a,b,c,d$ are elements of an $\SU$-matrix
$\left(\begin{array}{ll}
                 a & b \\
                 c & d 
              \end{array}\right)$
commuting with $z,\zb,dz,d\zb,\del,\delb$.
The above is a simpler notation of the equivalent one for left-coaction,
e.g., $\Delta_{L}(z) = (a\otimes z+b\otimes 1)(c\otimes z+d\otimes 1)^{-1}$.

The integration on $\S$, which is denoted by $<\cdot>_{\S}$,
can be defined as the integration on $\SU$ by
restricting the integrand to be an element in $\S$ and re-expressing it
through (\ref{S-SU}).

It is justified to call $dz d\zb (1+\zb z)^{-2}$ the volume form
in the sense that
if one treats the same whole algebra as the algebra on a quantum plane
so that everything else remains unchanged except $\del^{*} = -q^{2}\delb$,
the translational invariant integration on this plane $\int dz d\zb$ is related
to the integration $<\cdot>_{\S}$ by: 
\be
   \int dz d\zb (1+\zb z)^{-2} f = <f>_{\S} \label{int}
\eq
(up to normalization).

Details of everything above in this subsection can be found in \cite{CHZ}.
One can check that all the requirements of a complex manifold are met.

\subsection{ $\S$ As a Complex Manifold }

In this section we will treat $\S$ as a quantum complex manifold.

\subsubsection{ Metric And Connection } \label{MC}

The notation itself suggests that we take
$\Omega^{+}(\S) = \S dz$ and $\Omega^{-}(\S) = \S d\zb$.
Let the only possible value for an index be $0$, that is,
$z^{0}:=z$ and $\zb^{\bo}:=\zb$.
Also denote $g := g^{\bo 0}$ and $g^{-1} := g_{0\bo}$.
To define the metric $g$ we note that for $\S$ the Laplacian,
the $\SU$-invariant derivation of order two, is
\be
   \nabla^{2} = -c(1+\zb z)^{2}\delb\del, \label{lap}
\eq
where $c$ is an arbitrary real number.
On the other hand, the Laplacian of $\S$ as a complex manifold,
the holomorphic-transformation independent derivation of order two,
is $\del^{*}g\del$.
Equating these two expressions one gets
\[
   g = q^{2}c(1+\zb z)^{2}. \label{gc}
\]
Because the factor $(1+\zb z)$ will appear frequently,
we shall denote it in the following by:
\[
   \r := 1+\zb z.
\]

Classically any two-dimensional complex manifold is also
a K\"{a}hler manifold
and one can locally find a K\"{a}hler potential.
Analogy can be made here.
Define the K\"{a}hler form to be $K = dzg^{-1}d\zb$.
(It plays a special role in the differential calculus \cite{CHZ}.)
Obviously $dK = 0$ for the same reason as in the classical case.
The K\"{a}hler potential $V$ defined by $\d\bd V = K$ therefore exists.
One can solve $V$ in term of the deformed log:
\bee
   V &=& q^{-4}c^{-1}\sum_{n=0}^{\infty}\log_{q^{-1}}(1-\r^{-1}) \n \\
     &=& q^{-4}c^{-1}\sum_{n=0}^{\infty}\frac{\r^{-(n+1)}}{[n+1]_{q^{-1}}}, \n
\eqq
where $[n]_{q} = \frac{q^{2n}-1}{q^{2}-1}$.
In fact, the K\"{a}hler form $dz g^{-1} d\zb$ is just 
the volume form (up to normalization).

Using (\ref{omg1}) and (\ref{Rg1}),
we can immediately find the connection form and the curvature two-form:
\bee
   &\om_{0}^{\s 0} = -q^{2}(1+q^{2})dz\r^{-1}\zb, \label{om00} \\
   &R_{0}^{\s 0} = q^{4}(1+q^{2})dz d\zb \r^{-2}. \label{R00}
\eqq
It is easy to see that the torsion is zero in this case.

We define the vielbeins on $\S$ as
\bee
   &e := \r^{-1}dz, &\eb := \r^{-1}d\zb. \label{viel}
\eqq
The commutation relations between $e, \eb$ and $z, \zb$
are simply classical:
\bee
   &ez = ze, &e\zb =\zb e, \n \\
   &\eb z = z\eb, &\eb\zb = \zb\eb. \n
\eqq

The Hodge-$*$ map satisfying (\ref{starstar}) is given by
\bee
   &(*e) = ie, &(*\eb) = -i\eb, \n\\
   &(*1) = ic'^{-1}e\eb, &(*e\eb) = -ic', \label{c'}
\eqq
where $c'$ is a constant,
and the action of Hodge-$*$ on any form follows (\ref{hs}).

Now we can define $\d := -*d*$ and $\nabla^{2} := -\frac{1}{2}(d+\d)^{2}$.
\footnote{Here we have a factor of $\frac{1}{2}$
because we want to identify $\nabla^{2}$ with (\ref{lap}),
which is only the holomorphic part in $-(d+\d)^{2}$.}
When acting on a function $f$
\bee
   \nabla^{2}f &=& \frac{1}{2}*d*df \n \\
               &=& -c'\r^{2}\delb\del f. \n
\eqq
Hence $c'$ should be identified with $c$ because of (\ref{lap}).
This identification further justifies our Hodge-$*$ structure.

The scalar curvature is found to be
\be
   \R = cq^{2}(1+q^{2}). \label{S2SC}
\eq
The Ricci tensor defined by (\ref{hRT}) is
\be
   \R_0^{\s 0} = cq^4 (1+q^2). \label{S2RT}
\eq

\subsection{ $\S$ As a Riemannian Manifold }

In this section we treat $\S$ as a (real) Riemannian manifold.
The difference between Riemannian and complex manifolds
is that the latter is not
invariant under general coordinate transformation.
When indices are contracted for complex manifolds
they are summed over only half 
(the holomorphic part) 
of the possible values for an ordinary Riemannian manifold.

Assuming that $g^{00}=g^{\bo\bo}$ in the Riemannian case
(since we had $g^{\bo 0}=g^{0\bo}$ in the complex case),
we get from (\ref{lap})
\[
   g^{00}=g^{\bo\bo}=\frac{c}{1+q^{-2}}\r^{2}.
\]
Note that since the normalization of $g$ is changed
from Sec.(\ref{gc}),
the parameter $c'$ used in (\ref{c'}) should be changed
accordingly.

The equation $\nabla g^{00} = 0$ is identical to the one solved earlier 
for complex $\S$
and we assign $\om_{0}^{\s 0}$ to be the same as (\ref{om00}).
Let
\[
   \om_{\bo}^{\s \bo} = -(1+q^{-2})d\zb\r^{-1}z,
\]
which is the complex manifold connection form $\om_{0}^{\s 0}$ 
had we labelled the coordinates the opposite way: $\zb^{0}:=z$, 
and $z^{0}:=\zb$.
It solves $\nabla g^{\bo\bo} = 0$.

Similarly, the curvature two-form $R_{0}^{\s 0}$ is given by (\ref{R00})
and $R_{\bo}^{\s \bo}$ is
\[
   R_{\bo}^{\s \bo} = \d \om_{\bo}^{\s \bo} = (1+q^{-2})d\zb dz\r^{-2}.
\]

A straightforward calculation shows:
\bee
   &\R_{0}^{\s 0} = \frac{cq^{4}}{4}(1+q^{2})^{2}, \n \\
   &\R_{\bo}^{\s \bo} = \frac{cq^{-2}}{4}(1+q^{2})^{2}, \n \\
   &\R = \frac{c}{4}(1+q^{2})^{3}. \n
\eqq
This is different from (\ref{S2SC}) and (\ref{S2RT})
by constant factors.
The reason is that the Riemannian structure of $\S$
as a Riemannian manifold is concerned with
the general transformation and that of $\S$
as a complex manifold is concerned only with
the holomorphic transformations.
Unlike the situation in the classical case,
without leaving the holomorphic description of the complex manifold $\S$
one will never be able to know its Riemannian structure.
The discrepancy is introduced by the non-commutativity of the algebra.

\subsubsection{ Distances on $\S$ } \label{S2dis}

In this section we consider the ``distance'' between ``points'' on $\S$.
As mentioned earlier, (\ref{dis}) can be used to define a number
called the {\em distance} between any two states.

Before finding the distance we display a representation of $\S$
which shows clearly the correspondence between states and points.
The basis of the Hilbert space is labelled as $|k, \th\>$
for $k = 0, 1, 2, \cdots, \infty$, $\th \in [0, 2\pi)$.
(It is an irreducible representation given in \cite{Pod} for $\S$
supplemented with an additional index $\th$.
It is also equivalent by Fourier transform to
an irreducible representation given in \cite{Wor1} for $\SU$.)
The algebra is represented in the following way:
\been
   &z|k, \th\> = e^{i\th}(q^{-2k}-1)^{1/2}|k-1, \th\>, \\
   &\zb|k, \th\> = e^{-i\th}(q^{-2(k+1)}-1)^{1/2}|k+1, \th\>.
\eqqn
So we have
\[
   \r|k, \th\> = q^{-2k}|k, \th\>.
\]
Roughly speaking,
$\th$ corresponds to the azimuthal angel on $\S$ and
$\frac{q^{2k}-1}{q^{2k}+1}$ corresponds to the cosine of the polar angle.

However, for what follows it is not necessary to specify
the representation.
The only thing we need is $Sp(\r)$, the spectrum of $\r$,
which follows the commutation relations
\[
   z\r = q^{-2} \r z, \quad \zb \r = q^2 \r \zb,
\]
and that $\r=1+\zb z\geq1$.
It is easy to see that $Sp(\r) = \{q^{-2k}; k=0,1,2,\cdots\}$.
Hence in the following $\th$ is interpreted as
the collection of all parameters except $k$,
which labels the eigenvalue of $\r$.

Now we consider the distance between the two states
$|k, \th\>$ and $|k', \th\>$.
In the classical limit, it is just the radius ($\frac{1}{2}$) times
the difference in their polar angles:
\[
   D(p, p') = |F(z(p),\zb(p))-F(z(p'),\zb(p'))|,
\]
where $F(z,\zb) = \frac{1}{2}\cos^{-1}\left(\frac{z\zb-1}{z\zb+1}\right) 
= \sin^{-1}(\r^{-1/2})$.
For convenience we shall suppress the index $\th$ of a state from now on.
It is fixed for all consideration below.

Given the distance function $F$,
we can always decompose $F$ as $F = f(\r)+h(z,\zb)$,
where $h = \sum_{n=1}^{\infty} f_{n}(\r)z^{n}+g_{n}(\r)\zb^{n}$.
Since $\<k|h|k\>=0$ for all $k$,
if $F$ gives the distance between states $|k\>$ and $|k'\>$ then $f$ does, too.
But we have to check that the magnitude of $df$ is not larger than $1$.
Note that
\[
   |dF|^{2} = |df|^{2}+|dh|^{2}+(CT),
\]
where the cross-terms are
\[
   (CT)=(\del f)^{*}g(\del h)+(\delb f)^{*}g(\delb h)
        +(\del h)^{*}g(\del f)+(\delb h)^{*}g(\delb f).
\]
Since $\<p|(CT)|p\>=0$,
\[
   \<p|\left|df\right|^{2}|p\>\leq\<p|\left|dF\right|^{2}|p\>\leq
   \left\| |dF|^{2}\right\|\leq 1, \quad \forall p,
\]
which implies that $\left\| |df|^{2}\right\|\leq 1$.
Hence we only have to consider functions of $\r$
for our purpose.

Assume that the distance function is $F(\r)$ with $|dF|^{2} = 1$.
Because $(\del \r^{n}) = q\lam^{-1}\r^{-1}((q^{2}\r)^{n}-\r^{n})\zb$,
where $\lam = q-q^{-1}$, we have
\be
   (\del F(\r)) = q\lam^{-1}\r^{-1}(F(q^{2}\r)-F(\r))\zb. \label{delF}
\eq
Similarly,
\be
   (\delb F(\r)) = -q\lam^{-1}\r^{-1}(F(q^{-2}\r)-F(\r))z. \label{delbF}
\eq
Therefore,
\been
   |dF|^{2} &=& (\del F)^{*}g^{00}(\del F)+(\delb F)^{*}g^{\bo\bo}(\delb F) \\
            &=& \frac{cq^{4}}{\lam^{2}(1+q^{-2})}(W(q^{-2}\r)+W(\r)),
\eqqn
where $W(\r) = (\r-1)(F(q^{2}\r)-F(\r))^{2}$.
$|dF|^{2} = 1$ implies that
$W(\r)$ can only be the constant 
$\frac{1}{2}c^{-1}q^{-4}\lam^{2}(1+q^{-2})$.
Hence,
\be
   F(q^{2}\r)-F(\r) = -\lam\left(\frac{2cq^{4}(\r-1)}{1+q^{-2}}\right)^{-1/2}
   \label{FF}
\eq
and so $F$ can be solved as a power series expansion:
\[
   F(\r) = -\left(\frac{1+q^{-2}}{2cq^{2}}\right)^{1/2}
           \sum_{n=0}^{\infty}\frac{(2n)!}{(2^{n}n!)^{2}[n+1/2]_{q^{-1}}}
           \r^{-n-1/2}.
\]
This is not the only solution of (\ref{FF}).
Any function $f(\r)$ satisfying $f(q^2\r)=f(\r)$
can be added to it and (\ref{FF}) still holds.
However, due to the structure of $Sp(\r)$,
such functions will not contribute to the distance
between $|k\>$ and $|k'\>$.
For $q=1$, $c=4$, this solution is the power series expansion of
$-\sin^{-1}(\r^{-1/2})$.

It remains to argue that that the assumption $|dF|^{2} = 1$ is correct.
Consider a function $f(\r)$ with $|df|^{2} < 1$.
Then $F':=f+\eps F$ has $|dF'|^{2} < 1$ for $|\eps|$ sufficiently small.
And an appropriate phase of $\eps$ can make $|F'(k)-F'(k')| > |f(k)-f(k')|$.
So any $f$ with $|df|^{2} < 1$ is not the distance function.

Therefore we have the lemma:
\begin{lemma}
   The distance between states $|m, \th\>$ and $|n, \th\>$, $m\geq n$,
   according to (\ref{dis}) is equal to $F(q^{-2m})-F(q^{-2n})$.
\end{lemma}

The distance between the north pole and the south pole on $\S$,
for example, can be expressed as
\[
   F(\infty)-F(1) = \left(\frac{1+q^{-2}}{2cq^{2}}\right)^{1/2}
                      \sum_{n=0}^{\infty}\frac{(2n)!}{(2^{n}n!)^{2}
                      [n+1/2]_{q^{-1}}},
\]
which is the deformed $\pi/2$.

The distance between any two points can be obtained,
by using the quantum group symmetry of $SU_q(2)$,
from the distance
between the north pole ($z=\infty$) and an arbitrary point,
which we have just obtained above.

Using the commutation relations of $z, \zb$, one can check that
\be
   \T:=\left(
   \begin{array}{cc}
      z\r^{-1/2} & -q\r^{-1/2} \\
      \r^{-1/2} & \r^{-1/2}\zb
   \end{array}\right)
\eq
is an $SU_q(2)$-matrix.
This matrix transforms the north pole $(z=\infty)$ to $z$.
It also transforms the point
\[  z'':=(\d z'-q^{-1}\b)(-q\g z'+\a)^{-1}  \]
to $z'$.
The quantum group symmetry tells us that the distance between
$z$ and $z'$ is the same as the distance between the north pole and $z''$,
which is a function of $z, z'$.
Therefore we have:

\begin{prop}
   The distance between $(z, \zb)$ and $(z', \zb')$ on $\S$ is
   $|F(\rho'')|$,
   where $\rho''=(1+z''\zb'')=(1+z\zb)(1+z'\zb')(z-z')^{-1}(\zb-\zb')^{-1}$.
\end{prop}

Note that as the coordinates of points on the same sphere
the commutation relation between $z$ and $z'$ should be that of the
standard braiding \cite{CHZ2},
\[
   zz'=q^2 z'z - q\lam z'^2,
\] 
which is covariant under simultaneous $\SU$ transformation
on $z$ and $z'$.
This implies that $z$ and $z''$ simply commute with each other.
(The braiding is also formally satisfied by $(\infty, z'')$.
Divide the braiding relation on both sides by $z$ we get
$
z'z^{-1}=q^2 z^{-1}z'-q\lam z^{-1}z'^2 z^{-1}
$
which is satisfied by $(z,z')=(\infty,z'')$ but not by $(z,z')=(z'',\infty)$.)

A state  $|s\>$ in the Hilbert space representation of the braided algebra
generated by $\{z, \zb, z', \zb'\}$ corresponds to two ``points'' on $\S$.
So the distance between them is $\<s|F(\r'')|s\>$.
This is a modification of A. Connes' formula (\ref{dis})
requested by the braiding.

\subsubsection{ Connection with Connes' Formulation } \label{Hilbert}

Here we make a connection with A. Connes' 
quantum Riemannian geometry \cite{Con} by re-formulating
the quantum sphere in a way as close to his as possible.

To do so we consider the Hilbert space realization of $\O(\S)$.
The Hilbert space representation presented here is composed of two parts.
The first part $\{|\psi\>\}$ is the Hilbert space representing the algebra 
generated by $z, \zb$.
An example is given in Sec.\ref{S2dis}.
Another example is the GNS construction using the integration $<\cdot>$.
The second part $V$ is a vector space of, say,
2-component column vectors representing
the differential forms.
The differential calculus can then be represented in terms of 
the representation $\pi$
of $\S$ as (for $v\in V$):
\been
   &\pi(dz)|\psi>\tens v = \sqrt{c(1+q^{2})}\pi(\r)|\psi>\tens \tau v, \\
   &\pi(d\zb)|\psi>\tens v = \sqrt{c(1+q^{2})}\pi(\r)|\psi>\tens \tau^{\dagger} v,
\eqqn
where $\tau := \left(
                 \begin{array}{ll}
                     0 & 1 \\
                     0 & 0
                 \end{array}
              \right)$,
$\tau^{\dagger} := \left(
                 \begin{array}{ll}
                     0 & 0 \\
                     1 & 0
                 \end{array}
                 \right)$, and they satisfy:
$q\tau\tau^{\dagger}+q^{-1}\tau^{\dagger}\tau = \I$ for
        $\I := \left(
                 \begin{array}{ll}
                     q & 0 \\
                     0 & q^{-1}
                 \end{array}
              \right)$.
The vielbeins $e$ and $\eb$ (\ref{viel}) are represented
by the $\gamma$-matrices $\tau$ and $\tau^{\dagger}$,
which satisfy a deformed Clifford algebra.
The column $v$ is used to specify the direction of a cotangent vector at 
a ``point'' on $\S$.

Let the Dirac operator be
\[
   \Dir := k\left(\begin{array}{cc}
                   i & \zb \\
                   -z & -i
                \end{array}\right),
\]
where $k = q\lam^{-1}\sqrt{c(1+q^{2})}$.
It is chosen such that $dz = [D, z]$ and $d\zb = [D, \zb]$.
The goal is that the exterior derivative is realized by $D$.

Since $D^2$ is not central,
$(d^2\a) = [D^2, \a]$ for a form $\a$ is non-zero.
The nilpotency is achieved by taking the quotient of the algebra
over the ideal called the auxiliary fields.
They are the differential forms $\{a [D^2, b] c; a, b, c \in \A\}$.
For our case the auxiliary fields are found to be
\[
   a\left(\begin{array}{cc}   
              q & 0 \\
              0 & q^{-1}
          \end{array}\right)
\]
for all functions $a$ in $\S$.

The calculus is $\Z_{2}$-graded by
\[
   \g := k^{-2}(dz d\zb-d\zb dz)\r^{-2}
      = \left(\begin{array}{cc}   
            1 & 0 \\
            0 & -1
        \end{array}\right).
\]

The $\SU$-invariant integration on two-forms can be defined by the trace:
\[
   \int \a := Tr(\g\a |\Dir|^{-2}),
\]
where $Tr$ is the trace over the extended Hilbert space $\{|\psi>\tens v\}$.
Although this formula resembles that of A. Connes for 2-dimensional
integration \cite{Con}, according to him the power of $|\Dir|^{-1}$
in the trace is determined by the spectrum of $\Dir$.
In the classical case one gets the classical dimension,
but we get zero in this particular case.
Therefore unlike Connes' integration this one does not have 
the cyclic property $\int \a\b = \int \b\a$. 
Nevertheless,
it can be directly checked that the integration satisfies 
the consistency condition
\be
    \int Aux = 0 \label{aux}
\eq
for auxiliary fields, and the Stoke's theorem
\be
   \int d\a = 0 \label{Stoke}
\eq
for any one-form $\a$.
Stoke's theorem can be used to derive recursion relations for 
the integration of two-forms.
Equations (\ref{aux}) and (\ref{Stoke}) vanish already on 
the trace over the $2\times 2$ matrices
and hence remain so for any representation $|\psi>$ of the algebra $\S$.
They determine up to normalization the integration of two-forms.
Hence it agrees with the integration introduced before (\ref{int}).

\section{ The Complex Quantum Projective Spaces $CP_q(N)$ }\label{CPN}

The results in Sec.\ref{q-sphere} can be generalized to
the quantum projective spaces $CP_q(N)$ \cite{CHZ3}.
($S_q^2\sim CP_q(1)$.)
In Sec.\ref{CH} we consider the Hodge $*$ map.
In Sec.\ref{RG-CP} we find the Riemannian structure on $CP_q(N)$.

\subsection{ The Construction of the Hodge $*$ Map }\label{CH}

A prerequisite of the Riemannian structure is the Hodge $*$ map.
In general, if there exists for a quantum complex manifold
a K\"{a}hler form $K=dz^a g_{a\bb} d\zb^{\bb}$ which
is real and central, a Hodge $*$ satisfying (\ref{hs}) and (\ref{starstar})
can always be constructed.
Let
\bee
   {*(dz^{a_1}\cdots dz^{a_p})} :=
     dz^{a_1}\cdots dz^{a_p} K^{N-p}, \n\\
   {*(d\zb^{\ab_r}\cdots d\zb^{\ab_1})} :=
     K^{N-r} d\zb^{\ab_r}\cdots d\zb^{\ab_1}. \n
\eqq
Since $K$ is central, the property (\ref{hs}) is satisfied.
Since $K$ is real, (\ref{starstar}) is also satisfied.
Now we consider the Hodge $*$ of a differential form
which is not purely holomorphic and antiholomorphic.
The idea is to ``patch'' the holomorphic part
and the antiholomorphic part together.

Denote
\[
   \xi_a := g_{a\bb}d\zb^{\bb}.
\]
Because $K = dz^a\xi_a$ is central,
\[
K^p = (dz^{a_1}\cdots dz^{a_p})(\xi_{a_p}\cdots\xi_{a_1}).
\]
So we have
\[
   *(dz^{a_1}\cdots dz^{a_p}) = (dz^{a_1}\cdots dz^{a_N})
    (\eps_{a_1\cdots a_N}\xi_{a_N}\cdots\xi_{a_{p+1}}),
\]
where $\eps_{a_1\cdots a_N}$ is defined by
\[
   dz^{a_1}\cdots dz^{a_p} = \eps_{a_1\cdots a_N}dz^{a_1}\cdots dz^{a_N}.
\]
Similarly,
\[
   *(d\zb^{\ab_r}\cdots d\zb^{\ab_1}) = (\eps_{a_1\cdots a_N}
     \eta_{\ab_{r+1}} \cdots\eta_{\ab_N})(d\zb^{\bar{N}}\cdots d\zb^{\bar{1}}),
\]
where
\[
   \eta_{\ab} := (\xi_{a})^* = dz^b g_{b\ab}.
\]
Let $\mu(z, \zb)$ be the real function defined by the volume form:
\be \label{volume}
   K^N = dz^1\cdots dz^N \mu(z, \zb) d\zb^{\bar{N}}\cdots d\zb^{\bar{1}}.
\eq
Then we define
\be \label{patch}
   *(d\zb^{\bb_1}\cdots d\zb^{\bb_r} dz^{a_1}\cdots dz^{a_p})
   := (\eps_{b_1\cdots b_N}\eta_{\bb_{r+1}} \cdots\eta_{\bb_N})\mu^{-1}
   (\eps_{a_1\cdots a_N}\xi_{a_N}\cdots\xi_{a_{p+1}}),
\eq
Roughly speaking, we put the Hodge $*$ of the antiholomorphic part
and that of the holomorphic part together, and then take out from
the middle the volume form (\ref{volume}).
It can be shown that the properties (\ref{hs}) and (\ref{starstar}) hold.

The commutativity of the Hodge $*$ with all functions (\ref{hs})
is not necessary for the invariance of the scalar curvature.
As mentioned in Sec.\ref{Complex},
simply a prescription of ordering $*(f(z)\a g(\zb)) = f(z)(*\a)g(\zb)$
is sufficient.
Its significance is that any other prescription of ordering,
say, $*(f(\zb)\a g(z)) = f(\zb)(*\a)g(z)$,
gives the same result.

This construction of the Hodge $*$ map is, however, not unique.
When one patches the holomorphic and antiholomorphic parts
as in (\ref{patch}), one can choose to put the holomorphic part
before or after the antiholomorphic part.
They are in general inequivalent.
There is also the freedom to normalize (\ref{patch})
by different constant factors for each pair of $(p,r)$.

\subsection{ The Riemannian Structure on $CP_q(N)$ }\label{RG-CP}

The algebra of $CP_q(N)$ \cite{CHZ3} is given by the commutation relations:
\bee
   &z^a z^b = q^{-1}\hat{R}^{ab}_{cd}z^c z^d, \n\\
   &\zb^{\ab} z^b = q^{-1}(\hat{R}^{-1})^{bd}_{ac}z^c \zb^{\bar{d}}
     -q^{-1}\lam\d^b_a, \n\\
   &z^a dz^b = q\hat{R}^{ab}_{cd}dz^c z^d, \n\\
   &\zb^{\ab}dz^b = q^{-1}(\hat{R}^{-1})^{bd}_{ac}dz^c\zb^{\bar{d}}, \n
\eqq
where $\hat{R}$ is the $\hat{R}$-matrix of $GL_q(N)$ \cite{FRT}.
The $*$-involution is $z^{a*} = \zb^{\ab}$.

The K\"{a}hler form $K=dz^a g_{a\bb} d\zb^{\bb}$ for $CP_q(N)$ \cite{CHZ3}
is given by the deformed Fubini-Study metric
\[
   g_{a\bb} = q^{-1}\r^{-2}(\r\d_{ab}-q^2 \zb^a z^b),
\]
where $\r = 1+\sum_{a=1}^{N} z^a\zb^{\ab}$.
The inverse of the metric is $g^{\ab b} = q\r(\d_{ab}+\zb^{\ab} z^b)$.
This K\"{a}hler form is not only real and central,
but also invariant under the quantum group transformation.
\[
   z^a \rightarrow (T^0_0+z^b T_b^0)^{-1}(T_0^a+z^c T_c^a),
\]
where $T_a^b$ is an $SU_q(N+1)$-matrix.
Consequently its corresponding Hodge $*$ defined as above
is also commutative with this quantum group transformation.

The deformed Fubini-Study metric implies that the connection one-form is
\[
   \om_a^{\s b} = C_{ac}^{bd} \zb^{\bar{c}}\r^{-1}dz^d,
\]
where $C_{ac}^{bd} = \d_{ac}\d_{bd} + q^{N-d}\d_{ab}\d_{cd}$.
The curvature two form is
\[
   R_a^{\s b} = -C_{ac}^{bd} g_{c\bar{e}}d\zb^{\bar{e}}dz^d.
\]

The scalar curvature and Ricci tensor are, up to normalization,
\be
   \R \propto 1, \quad \R_a^{\s b} \propto \d_a^b. \n
\eq
As in the classical case,
this result can also be obtained by arguments
based on the quantum group symmetry.

\section{ The Two-Sheeted Space }\label{Sheets}

Using the algebraic formulation of Riemannian geometry,
we reproduce in this section the theory of gravity
for the two-sheeted space which was first described in Ref.\cite{CFF}.
In that paper the Riemannian geometry on quantum spaces
is formulated in terms of A. Connes' non-commutative geometry.

The two-sheeted space is the product of a classical 4-dimensional
manifold ${\cal M}_4$ and a space of two discrete points $\Z_2$.
Denote the two points in $\Z_2$ as $a$ and $b$.
The algebra of functions on $\Z_2$ is generated by $1$ and $e$,
where $1(a)=1(b)=1$ and $e(a)=-e(b)=1$.
It follows that $e$ is real and
\be
   e^2 = 1. \label{e2}
\eq
Let the exterior derivative to act on (\ref{e2}) we find
\[
   e de = -de e.
\]
We also define $dede=0$.
In addition, $e$ commutes with the coordinates $\{x^{\mu}\}$
on ${\cal M}_4$ and $\{dx^{\mu}\}$,
and $de$ commutes with $\{x^{\mu}\}$ and anti-commutes with
$\{dx^{\mu}\}$.

To obtain the results in Ref.\cite{CFF} we assume that the vielbeins
and connection one-forms can be written as
\[
   E^a=dx^{\mu}e_{\mu}^a, \quad E^5=de\lam \n
\]
and
\[
   \Omega^{AB}=dx^{\mu}(\om^{AB}_{\mu}+ev^{AB}_{\mu})+de(l^{AB}+ek^{AB}),
\]
where $e_{\mu}^a$, $\lam$, $\om^{AB}_{\mu}$, $v^{AB}_{\mu}$, $l^{AB}$ and $k^{AB}$
are all real functions of $x$.
The indices $A, B$ take values in $\{1,2,3,4,5\}$,
where $\{1,2,3,4\}$ correspond to $dx^{\mu}$ or $E^a$,
and $\{5\}$ corresponds to $de$ or $E^5$.

The Hodge $*$ map defined on $E^A$ is the classical one.
For example, $*(E^A E^B)=\frac{1}{3!}\epsilon^{ABCDEF}E^C E^D E^E E^F$.
This map does not have the property (\ref{hs}),
hence the Lagrangian (\ref{Lag}) is invariant only under
the coordinate transformation restricted to ${\cal M}_4$,
i.e., $x^{\mu}\rightarrow x'^{\mu}(x), e\rightarrow e$.

The integration over the whole space can be decomposed into the usual integration
over the four-dimensional manifold followed by the integration
\[
   \int_{Z_2}E^5(a+be)=a
\]
for arbitrary numbers $a,b$.
(The requirement that $\int_{Z_2}dee$ vanishes
implies the cyclic property of Connes' integration in this case:
$\int\a\b=\int\b\a$.)

Using the metricity condition (\ref{nablag}) and
the torsion-free condition (\ref{tor}) one can partially
solve for the connection.
But many components of the connection are still free.
They should be viewed as independent fields.
It turns out that they are not dynamical fields
because in the Lagrangian (the scalar curvature) they
do not have time derivatives.
Their equations of motion are simply constraints which
are solved by their vanishing.

The action for the gravity on this two-sheeted space defined by
\[
   I = \int_{{\cal M}_4\times Z_2} E^A(*R^{AB})E^B
\]
is, after taking out all non-dynamical fields,
the same as \cite{CFF}
\[
   I = - \int_{{\cal M}_4}(\R_4 - 2\lam^{-1}\nabla_{\mu}\del^{\mu}\lam)\sqrt{g_4}d^4 x,
\]
where $\R_4$ is the usual scalar curvature of ${\cal M}_4$,
$\nabla_{\mu}$ is the usual covariant derivative on ${\cal M}_4$
and $g_4$ is the determinant of the metric on ${\cal M}_4$.
The only new dynamical field introduced by the $\Z_2$ structure in spacetime is $\lam(x)$.
By changing variable $\lam=exp(\sigma)$ \cite{CFF} we get
\[
   I = - \int_{{\cal M}_4}(\R_4 - 2\del_{\mu}\sigma\del^{\mu}\sigma)\sqrt{g_4}d^4 x.
\]

\section{ Conclusion }

In this paper we proposed a straightforward formulation of 
Riemannian geometry
on quantum spaces with a $*$-involution and a Hodge $*$ map,
and we showed several examples.
In addition to the possibility of applying it
to describe physics at the Planck scale,
this formulation can be used for Kaluza-Klein theories
to build models with the extra dimensions corresponding to quantum spaces.

\section{ Acknowledgements }

The author would like to thank Professor Bruno Zumino for invaluable 
advices, discussions, encouragement and support.
The author also appreciate the discussions with Piotr Podle\'{s}
and Chong-Sun Chu.

This work was supported in part by the Director, Office of
Energy Research, Office of High Energy and Nuclear Physics, Division of
High Energy Physics of the U.S. Department of Energy under Contract
DE-AC03-76SF00098 and in part by the National Science Foundation under
grant PHY-90-21139.

\end{document}